\documentclass[11pt]{amsart}
\usepackage{amsmath,amssymb,mathrsfs}
\newtheorem{theorem}{Theorem}

\newtheorem{corollary}[theorem]{Corollary}
\newtheorem{lemma}[theorem]{Lemma}

\begin{document}

\title{A local noncollapsing estimate for mean curvature flow}
\author{Simon Brendle and Keaton Naff}
\address{Department of Mathematics \\ Columbia University \\ New York NY 10027}
\address{Department of Mathematics \\ Columbia University \\ New York NY 10027}
\begin{abstract}
We prove a local version of the noncollapsing estimate for mean curvature flow. By combining our result with earlier work of X.-J. Wang, it follows that certain ancient convex solutions that sweep out the entire space are noncollapsed. 
\end{abstract}
\thanks{The first author was supported by the National Science Foundation under grant DMS-1806190 and by the Simons Foundation.}
\maketitle

\section{Introduction}

A crucial property of mean curvature flow is that embeddedness is preserved under the evolution, and this property is especially consequential in the class of mean-convex flows. Indeed, White \cite{Whi00,Whi03} proved several important results on singularities of embedded, mean-convex flows. Among other things, White showed certain ``collapsed'' singularity models, like the grim reaper or the multiplicity-two hyperplane, cannot arise as blow-up limits of embedded mean-convex flows.

In \cite{SW09}, Sheng and Wang introduced a quantitative version of the concept of embeddedness. Let $M_t$ be a family of embedded, mean-convex hypersurfaces evolving by mean curvature flow. We say that the flow $M_t$ is $\alpha$-noncollapsed if 
\[\frac{1}{2} \, \alpha^{-1} \, H(x,t) \, |x-y|^2 - \langle x-y,\nu(x,t) \rangle \geq 0\] 
for all points $x,y \in M_t$. This concept has a natural geometric interpretation. A flow is $\alpha$-noncollapsed if, for every point $x \in M_t$, there exists a ball of radius $\alpha H(x,t)^{-1}$ which lies in the inside of $M_t$ and which touches $M_t$ at $x$. 

It is a consequence of White's work that every compact, embedded, mean-convex solution of mean curvature flow is $\alpha$-noncollapsed, for some $\alpha > 0$, up to the first singular time. Alternative proofs of this result were given by Sheng and Wang \cite{SW09} and by Andrews \cite{And12}. In \cite{Bre15}, the first author proved a sharp version of this noncollapsing estimate. More precisely, if we start from a closed, embedded, mean-convex solution of mean curvature flow, then every blow-up limit is $1$-noncollapsed.

In this note, we prove a local version of the noncollapsing estimate for the mean curvature flow. 

\begin{theorem}
\label{main.thm}
Let us fix radii $R$ and $r$ such that $R \geq \sqrt{1+3n} \, r$. Moreover, let $\Lambda$ be a positive real number. Let $M_t$, $t \in [-r^2,0]$, be an embedded solution of mean curvature flow in the ball $B_{3R}(0)$ satisfying $R^{-1} \leq \Lambda H$ and $|A| \leq \Lambda H$ for all $t \in [-r^2,0]$ and all $x \in M_t \cap B_{\sqrt{1+3n} \, r}(0)$. Moreover, suppose that 
\[\frac{1}{2} \, \Lambda \, H(x,-r^2) \, |x-y|^2 - \langle x-y,\nu(x,-r^2) \rangle \geq 0\] 
for all $x \in M_{-r^2} \cap B_{\sqrt{1+3n} \, r}(0)$ and all $y \in M_{-r^2} \cap B_{3R}(0)$. Then
\[\Lambda \, (2+6n)^{4} \, H(x,t) \, |x-y|^2 - \langle x-y,\nu(x,t) \rangle \geq 0\]  
for all $t \in [-r^2,0]$, all $x \in M_t \cap B_{\frac{r}{2}}(0)$, and all $y \in M_t \cap B_{3R}(0)$. 
\end{theorem}

The proof of Theorem \ref{main.thm} relies on two main ingredients. The first is the evolution equation from \cite{Bre15} for the reciprocal of the inscribed radius. The second is a particular choice of cutoff function inspired in part by the work of Ecker and Huisken \cite{EH91}. The argument can be readily adapted to fully nonlinear flows given by speeds $G = G(h_{ij})> 0$ which are homogeneous of degree one, concave, and satisfy $0 < \frac{\partial G}{\partial h_{ij}} \leq K g_{ij}$ for a uniform constant $K$. 

Let us now discuss an application of the main theorem. 

\begin{corollary}
\label{ancient.solution.1}
Let $\Lambda$ be a positive real number. Let $M_t$, $t \in (-\infty,0]$, be an embedded ancient solution of mean curvature flow in $\mathbb{R}^{n+1}$ such that $H > 0$ and $|A| \leq \Lambda H$ at each point in space-time. Suppose that there exists a sequence of times $t_j \to -\infty$ such that 
\[\frac{1}{2} \, \Lambda \, H(x,t_j) \, |x-y|^2 - \langle x-y,\nu(x,t_j) \rangle \geq 0\] 
for all $x \in M_{t_j} \cap B_{\sqrt{-(1+3n)t_j}}(0)$ and all $y \in M_{t_j}$. Then 
\[\Lambda \, (2+6n)^{4} \, H(x,t) \, |x-y|^2 - \langle x-y,\nu(x,t) \rangle \geq 0\]
for all $t \in (-\infty,0]$, all $x \in M_t$, and all $y \in M_t$. 
\end{corollary}

To deduce Corollary \ref{ancient.solution.1} from Theorem \ref{main.thm}, we put $r_j = \sqrt{-t_j}$. Moreover, for each $j$, we choose $R_j$ large enough so that $R_j \geq \sqrt{-(1+3n)t_j}$ and $R_j^{-1} \leq \Lambda H$ for all $t \in [t_j,0]$ and all $x \in M_t \cap B_{\sqrt{-(1+3n)t_j}}(0)$. If we apply Theorem \ref{main.thm} and take the limit as $j \to \infty$, the assertion follows. 

\begin{corollary}
\label{ancient.solution.2}
Let $M_t$, $t \in (-\infty,0]$, be a convex ancient solution of mean curvature flow in $\mathbb{R}^{n+1}$ with $H > 0$. Suppose that there exists a sequence $t_j \to -\infty$ such that the rescaled hypersurfaces $(-t_j)^{-\frac{1}{2}} \, M_{t_j}$ converge in $C_{\text{\rm loc}}^\infty$ to a cylinder $S^{n-k} \times \mathbb{R}^k$ with multiplicity $1$, where $k \in \{0,1,\hdots,n-1\}$. Then there exists a constant $\Lambda(n)$ such that 
\[\Lambda(n) \, H(x,t) \, |x-y|^2 - \langle x-y,\nu(x,t) \rangle \geq 0\] 
for all $t \in (-\infty,0]$, all $x \in M_t$, and all $y \in M_t$. 
\end{corollary}

In \cite{Wan11}, X.-J. Wang considered ancient solutions to mean curvature flow that can be expressed as level sets $M_t = \{u=-t\}$, where $u$ is a convex function which is defined on the entire space $\mathbb{R}^{n+1}$ and satisfies 
\[\sum_{i,j=1}^{n+1} \Big ( \delta_{ij} - \frac{D_i u \, D_j u}{|\nabla u|^2} \Big ) \, D_i D_j u = 1.\] 
Wang showed that such ancient solutions admit a blow-down limit which is a cylinder $S^{n-k} \times \mathbb{R}^k$ with multiplicity $1$, where $k \in \{0,1,\hdots,n-1\}$ (see \cite{Wan11}, Theorem 1.3). By Corollary \ref{ancient.solution.2}, such an ancient solution must be noncollapsed. 

\section{Proof of Theorem \ref{main.thm}}

By scaling, it suffices to prove the assertion for $r=1$. Let us fix a radius $R \geq \sqrt{1+3n}$. Moreover, we fix a positive real number $\Lambda$. Let $M_t$, $t \in [-1,0]$, be an embedded solution of mean curvature flow in the ball $B_{3R}(0)$ satisfying $R^{-1} \leq \Lambda H$ and $|A| \leq \Lambda H$ for all $t \in [-1,0]$ and all $x \in M_t \cap B_{\sqrt{1+3n}}(0)$. Moreover, we assume that 
\[\frac{1}{2} \, \Lambda \, H(x,-1) \, |x-y|^2 - \langle x-y,\nu(x,-1) \rangle \geq 0\] 
for all $x \in M_{-1} \cap B_{\sqrt{1+3n} }(0)$ and all $y \in M_{-1} \cap B_{3R}(0)$. 

We define a function $\varphi$ by 
\[\varphi(x,t) := (2\Lambda)^{-\frac{1}{4}} \, (1+3n)^{-1} \, (1-|x|^2-3nt)\] 
for all $t \in [-1,0]$ and all $x \in M_t \cap B_{\sqrt{1-3nt}}(0)$. Moreover, we define 
\[\Phi(x,t) := \varphi(x,t)^{-4} \, H(x,t)\] 
for all $t \in [-1,0]$ and all $x \in M_t \cap B_{\sqrt{1-3nt}}(0)$. In the following,  let $\lambda_1 \leq \hdots \leq \lambda_n$ denote the eigenvalues of the second fundamental form.

\begin{lemma} 
\label{aux}
We have $\Phi \geq 2\Lambda H$ for all $t \in [-1,0]$ and all $x \in M_t \cap B_{\sqrt{1-3nt}}(0)$. Moreover, the eigenvalues of the second fundamental form satisfy $|\lambda_i| \leq \frac{\Phi}{2}$ for all $t \in [-1,0]$ and all $x \in M_t \cap B_{\sqrt{1-3nt}}(0)$.
\end{lemma} 

\textbf{Proof.} 
By definition, $\varphi \leq (2\Lambda)^{-\frac{1}{4}}$ and $\Phi \geq 2\Lambda H$ for all $t \in [-1,0]$ and all $x \in M_t \cap B_{\sqrt{1-3nt}}(0)$. This proves the first statement.

To prove the second statement, we observe that $|A| \leq \Lambda H \leq \frac{\Phi}{2}$ for all $t \in [-1,0]$ and all $x \in M_t \cap B_{\sqrt{1-3nt}}(0)$. This completes the proof of Lemma \ref{aux}. \\

\begin{lemma}
\label{testfunction}
The function $\varphi$ satisfies 
\[\frac{\partial \varphi}{\partial t} - \Delta \varphi < 0\] 
for all $t \in [-1,0]$ and all $x \in M_t \cap B_{\sqrt{1-3nt}}(0)$.
\end{lemma}

\textbf{Proof.} 
We compute 
\[\frac{\partial \varphi}{\partial t} - \Delta \varphi = -(2\Lambda)^{-\frac{1}{4}} \, (1+3n)^{-1} \, n < 0 \] 
for all $t \in [-1,0]$ and all $x \in M_t \cap B_{\sqrt{1-3nt}}(0)$. \\

\begin{lemma} 
\label{inequality.for.Phi}
The function $\Phi$ satisfies
\[\frac{\partial \Phi}{\partial t} - \Delta \Phi - |A|^2 \, \Phi + 2 \sum_{i=1}^n \frac{(D_i \Phi)^2}{\Phi-\lambda_i} > 0\] 
for all $t \in [-1,0]$ and all $x \in M_t \cap B_{\sqrt{1-3nt}}(0)$.
\end{lemma}

\textbf{Proof.}
Using Lemma \ref{testfunction}, we obtain 
\begin{align*} 
\frac{\partial \Phi}{\partial t} - \Delta \Phi
&= \varphi^{-4} \, \Big ( \frac{\partial H}{\partial t} - \Delta H \Big ) - 4 \, \varphi^{-5} H \, \Big ( \frac{\partial \varphi}{\partial t} - \Delta \varphi \Big ) \\ 
&- \frac{4}{3} \, \varphi^{-4} H \, \Big | \frac{\nabla H}{H} - 4 \, \frac{\nabla \varphi}{\varphi} \Big |^2 + \frac{4}{3} \, \varphi^{-4} H \, \Big | \frac{\nabla H}{H} - \frac{\nabla \varphi}{\varphi} \Big |^2 \\ 
&> |A|^2 \, \varphi^{-4} H - \frac{4}{3} \, \varphi^{-4} H \, \Big | \frac{\nabla H}{H} - 4 \, \frac{\nabla \varphi}{\varphi} \Big |^2 \\ 
&= |A|^2 \, \Phi - \frac{4}{3} \, \frac{|\nabla \Phi|^2}{\Phi} 
\end{align*}
for all $t \in [-1,0]$ and all $x \in M_t \cap B_{\sqrt{1-3nt}}(0)$. Moreover, it follows from Lemma \ref{aux} that $\frac{\Phi}{2} \leq \Phi-\lambda_i \leq \frac{3\Phi}{2}$ for all $t \in [-1,0]$ and all $x \in M_t \cap B_{\sqrt{1-3nt}}(0)$. Consequently, 
\[\frac{\partial \Phi}{\partial t} - \Delta \Phi > |A|^2 \, \Phi - \frac{4}{3} \, \frac{|\nabla \Phi|^2}{\Phi} \geq |A|^2 \, \Phi - 2 \sum_{i=1}^n \frac{(D_i \Phi)^2}{\Phi-\lambda_i}\] 
for all $t \in [-1,0]$ and all $x \in M_t \cap B_{\sqrt{1-3nt}}(0)$. This completes the proof of Lemma \ref{inequality.for.Phi}. \\

We next define 
\[Z(x,y,t) := \frac{1}{2} \, \Phi(x,t) \, |x-y|^2 - \langle x-y,\nu(x,t) \rangle\] 
for $t \in [-1,0]$, $x \in M_t$, and $y \in M_t$. \\

\begin{lemma}
\label{two.point.function}
We have $Z(x,y,t) \geq 0$ for all $t \in [-1,0]$, all $x \in M_t \cap B_{\sqrt{1-3nt}}(0)$, and all $y \in M_t \cap B_{3R}(0)$. 
\end{lemma}

\textbf{Proof.} 
Suppose that the assertion is false. Let $J$ denote the set of all times $t \in [-1,0]$ with the property that we can find a point $x \in M_t \cap B_{\sqrt{1-3nt}}(0)$ and a point $y \in M_t \cap B_{3R}(0)$ such that $Z(x,y,t) < 0$. Moreover, we define $\bar{t} := \inf J$.

By definition, $Z(x,y,t) \geq 0$ for all $t \in [-1,\bar{t})$, all $x \in M_t \cap B_{\sqrt{1-3nt}}(0)$, and all $y \in M_t \cap B_{3R}(0)$. Moreover, we can find a sequence of times $t_j \in J$ such that $t_j \to \bar{t}$. For each $j$, we can find a point $x_j \in M_{t_j} \cap B_{\sqrt{1-3nt_j}}(0)$ and a point $y_j \in M_{t_j} \cap B_{3R}(0)$ such that $Z(x_j,y_j,t_j) < 0$. Clearly, $x_j \neq y_j$, and 
\[\frac{2 \, \langle x_j-y_j, \nu(x_j,t_j) \rangle}{|x_j -y_j|^2} > \Phi(x_j,t_j).\] 
Using the Cauchy-Schwarz inequality, we obtain 
\[\frac{2}{|x_j-y_j|} > \Phi(x_j,t_j) \geq 2\Lambda H(x_j,t_j) \geq 2R^{-1}.\] 
In other words, $|x_j-y_j| \leq R$. After passing to a subsequence, the points $x_j$ converge to a point $\bar{x} \in M_{\bar{t}}$ satisfying $|\bar{x}| \leq \sqrt{1-3n\bar{t}}$. Moreover, the points $y_j$ converge to a point $\bar{y} \in M_{\bar{t}}$ satisfying  $|\bar{x}-\bar{y}| \leq R$. In particular, $|\bar{y}| \leq \sqrt{1+3n} + R \leq 2R$. 

If $|\bar{x}| = \sqrt{1-3n\bar{t}}$ and $\bar{x} \neq \bar{y}$, then 
\[\frac{2 \, \langle \bar{x}-\bar{y}, \nu(\bar{x},\bar{t}) \rangle}{|\bar{x}-\bar{y}|^2} = \limsup_{j \to \infty} \frac{2 \, \langle x_j-y_j, \nu(x_j,t_j) \rangle}{|x_j -y_j|^2} \geq \limsup_{j \to \infty} \Phi(x_j,t_j) = \infty,\] 
which is impossible. 

If $|\bar{x}| = \sqrt{1-3n\bar{t}}$ and $\bar{x} = \bar{y}$, then 
\[\lambda_n(\bar{x},\bar{t}) \geq \limsup_{j \to \infty} \frac{2 \, \langle x_j-y_j, \nu(x_j,t_j) \rangle}{|x_j -y_j|^2} \geq \limsup_{j \to \infty} \Phi(x_j,t_j) = \infty,\] 
which is impossible. 

If $|\bar{x}| < \sqrt{1-3n\bar{t}}$ and $\bar{x} = \bar{y}$, then 
\[\lambda_n(\bar{x},\bar{t}) \geq \limsup_{j \to \infty} \frac{2 \, \langle x_j-y_j, \nu(x_j,t_j) \rangle}{|x_j -y_j|^2} \geq \limsup_{j \to \infty} \Phi(x_j,t_j) = \Phi(\bar{x},\bar{t}),\] 
which contradicts Lemma \ref{aux}. 

Therefore, we must have $|\bar{x}| < \sqrt{1-3n\bar{t}}$ and $\bar{x} \neq \bar{y}$. Moreover,  
\[\frac{2 \, \langle \bar{x}-\bar{y}, \nu(\bar{x},\bar{t}) \rangle}{|\bar{x}-\bar{y}|^2} = \limsup_{j \to \infty} \frac{2 \, \langle x_j-y_j, \nu(x_j,t_j) \rangle}{|x_j -y_j|^2} \geq \limsup_{j \to \infty} \Phi(x_j,t_j) \geq \Phi(\bar{x},\bar{t}).\] 
We claim that $\bar{t} \in (-1,0]$. Indeed, if $\bar{t}=-1$, then our assumption implies 
\[\frac{2 \, \langle \bar{x}-\bar{y}, \nu(\bar{x},\bar{t}) \rangle}{|\bar{x}-\bar{y}|^2} \leq \Lambda \, H(\bar{x},\bar{t}) < \Phi(\bar{x},\bar{t}),\] 
which is impossible. Consequently, $\bar{t} \in (-1,0]$.

To summarize, we have shown that $\bar{t} \in (-1,0]$, $Z(\bar{x},\bar{y},\bar{t}) \leq 0$, and $Z(x,y,t) \geq 0$ for all $t \in [-1,\bar{t})$, all $x \in M_t \cap B_{\sqrt{1-3nt}}(0)$, and all $y \in M_t \cap B_{3R}(0)$. Arguing as in the proof of Proposition 2.3 in \cite{Bre15}, we conclude that 
\[\frac{\partial \Phi}{\partial t} - \Delta \Phi - |A|^2 \, \Phi + 2 \sum_{i=1}^n \frac{(D_i \Phi)^2}{\Phi - \lambda_i} \leq 0\] 
at the point $(\bar{x},\bar{t})$. This contradicts Lemma \ref{inequality.for.Phi}. This completes the proof of Lemma \ref{two.point.function}. \\

We now complete the proof of Theorem \ref{main.thm}. By definition, $\varphi(x,t) \geq (2\Lambda)^{-\frac{1}{4}} \, (2+6n)^{-1}$ and $\Phi(x,t) \leq 2\Lambda \, (2+6n)^{4} \, H(x,t)$ for all $t \in [-1,0]$ and all $x \in M_t \cap B_{\frac{1}{2}}(0)$. Using Lemma \ref{two.point.function}, we conclude that 
\[\Lambda \, (2+6n)^{4} \, H(x,t) \, |x-y|^2 - \langle x-y,\nu(x,t) \rangle \geq Z(x,y,t) \geq 0\] 
for all $t \in [-1,0]$, all $x \in M_t \cap B_{\frac{1}{2}}(0)$, and all $y \in M_t \cap B_{3R}(0)$. This completes the proof of Theorem \ref{main.thm}.

\end{document}